\documentclass{llncs}
\usepackage{makeidx}  
\usepackage[pdftex]{graphicx}
\usepackage{float}
\usepackage{amsmath}
\usepackage {comment}
\usepackage {framed}
\usepackage {xcolor}
\colorlet{shadecolor}{yellow!50}
\usepackage{ulem}
\begin{document}
\frontmatter          
\pagestyle{headings}  
\addtocmark{Hamiltonian Mechanics} 
\mainmatter              
\title{ Development of the method of computer analogy for studying and solving complex nonlinear systems}
\titlerunning{Computer Analogy}  
%
\author{Vladimir Aristov \and Andrey Stroganov }
\authorrunning{Vladimir Aristov, Andrey Stroganov} 
%
\tocauthor{Vladimir Aristov, Andrey Stroganov}
\institute{Dorodnicyn Computing Centre of Russian Academy of Sciences,\\
\email{aristovvl@yandex.ru},
\email{savthe@gmail.com}
}

\maketitle              

\begin{abstract}
A method of representation of a solution as segments of the series in powers of the step of the 
independent variable is expanded for solving complex systems of ordinary differential equations (ODE): 
the Lorenz system and other systems. A new procedure of reduction of the representation of the solution 
to a sum of two parts (regular and random)  is performed.  A shifting procedure is applied in each level 
of the independent variable
to the random part and it acts as the filter that extracts the values to the regular part. 
In certain cases it is possible 
to omit the random part and construct the approximation which does not converge but still provides the 
qualitative information about the full solution (a linear approximation provides a simple exact 
solution). Evaluation of the error for this case is performed. Constructing the analytical representation 
of the solutions for these systems by the developed method is presented. 

\keywords{theoretical computer model; digit shifting;  
the method of computer analogy; nonlinear differential equations }
\end{abstract}

\section{Introduction}
In papers [1-2] we proposed a new approach for solving nonlinear differential 
equations which utilizes two main principles of storing the numbers in the digital computers: 
1) the numbers are represented as the segments of a power series; 2) there is a digit shifting procedure. 
In each level the numerical solution is presented as the 
segment of a series in the powers of the step $\tau$ of the independent variable. 
The digit shifting procedure is applied in each level to change the values of the coefficients
of the terms so each coefficient would not exceed $\tau^{-1}$. We have shown that the
digit shifting can produce quasi-random numbers which have been used 
to exclude the intermediate levels of computations.
This method can provide a solution in the explicit form (as a computer provides a solution in the 
numerical form, after executing many intermediate and "hidden" operations). 
Some examples of solving different nonlinear equations and simple systems have been demonstrated. 

In the present paper we consider nonlinear systems which have the physical sense. 
Namely the Van der Pol equation (equivalent to a system of two ODE) and the well-known Lorenz system 
of the three ODE. We prove that for ODE where the right hand part is a polynomial of degree $m$, it is 
sufficient to retain the terms up to $m$ plus the order of accuracy of the numerical method.
	
\section{The method of computer analogy}
\subsection {Problem statement}
In most cases the nonlinear differential equations can not be solved analytically, but usually they
can be solved numerically instead. However one of the main distinctions between analytical and numerical 
solutions is that the numerical solution is not conceivable without a computer --- a physical device 
which performs routine computations. 

Consider the following Cauchy problem:
\begin{equation}\label{eq:ode}
\frac{dy}{dt} = f(y), \quad y(0) = A.
\end{equation}

Let there be a $q$-th order explicit finite-difference scheme which approximates the solution, 
it can be written as:
\[
y_{n+1} = y_n + \tau G(y_n) + O(\tau^{q+1}), \quad y_0 = A,
\]
where $\tau$ is a step of the independent variable $t$, $G(y)$ is determined by function a $f$ and by the
chosen finite difference method. Note that for the first order method $G$ is equivalent to $f$.

When solving the numerical solution the computer acts as a black box, which  
runs an algorithm and after performing a large amount of steps of the algorithm provides 
a result as the array of numbers. The computer operations are fairly 
simple, but due to nonlinearity (in general) of the problem under consideration 
it is difficult to understand how the values 
will be changed after few operations. 

From a certain point of view the computer solves the numerical analog of the analytical problem, and 
instead of the analysis and analytical computations it executes an algorithm. But on 
the other hand we can consider a numerical solution as the independent problem and 
search for an analytical solution (or approximation) for this particular problem. 
First, we will formalize the numbers representation in computers. 
The following
aspects are essential: 1) the numbers are stored as the segments of the power series, 2) there is
the digit shifting procedure which shifts a part of the value to the left digit.

\subsection {Solution representation}

The power series in powers of $\tau$ can be obtained from the explicit finite differece scheme and the 
Taylor expansion. 
\begin{theorem}
Let $G(y)$ is infinitely differentiable function with bounded derivatives. Then $y_n$ can be 
represented by the series in powers of $\tau$:
\begin{equation}\label{eq:y_series}
y_n = \sum_{m = 0}^{\infty} a_{m,n}\tau^m.
\end{equation}
\end{theorem}
This follows from the Taylor expansion of $G(y_n)$ at $y^*=a_{0,n}$:
\[
G(y_n) = G(y^*) + \sum_{m=1}^{\infty}\frac{G^{(m)}(y^*)}{m!}\left(y_n-y^*\right)^m.
\]

Using this expansion we find $y_{n+1}$. Comparison of the coefficients in the 
sequential layers gives the following recurrent relations:
\begin{equation}\label{eq:recurrent_coefficients}
\begin{array}{l}
a_{0,n+1} = a_{0,n}, \\
a_{1,n+1} = a_{1,n} + G(y^*), \\
a_{2,n+1} = a_{2,n} + a_{1,n}G'(y^*), \\
a_{3,n+1} = a_{2,n} + a_{2,n}G'(y^*) + \frac{1}{2}a_{1,n}^2G''(y^*),
\end{array}
\end{equation}
and so on.

The expansion of $y_n$ can not be truncated after a fixed number of terms because the 
truncated value will be of a greater order than the accuracy of the finite-difference method, and such
a solution would diverge.

\begin {comment}
Let us see the two basic principles of storing the numbers.

1. The numbers are stored as the segments of the power series. For classical computer those are the
powers of 2. Cosider the number $11$ in binary numeral system. It is $1011_2$, where each digit 
denotes the quantity of powers of $2$. The position of digit is also important, 
it defines the qualitative
''weight'' of the digit, so the left unity corresponds to $2^3$, the zero corresponds 
to $2^2$ and so on. 
Number $11$ can be represented in binary system:
\[
1011_2 = 1\cdot 2^3 + 0\cdot 2^2 + 1 \cdot 2^1 + 1 \cdot 2^0 = 11_{10}
\]

2. if there is a potential overflow in some digit, the part of it is shifted to the left digit (which 
corresponds to the value of the higher order). Let us add to 11 the number 2. In binary numeral system this
implies $1011 + 10$.
\[
11_{10} + 2_{10} = 1011_2 + 10_2 = 1101_2
\]

We will formalize these principles and apply them for solving nonlinear differential equations.

\subsection{Problem definition}

\subsection {The expansion of the solution in powers of $\tau$}
Notice that the finite difference scheme can can generate series in powers of $\tau$.
\begin{theorem}
Let $G(y)$ is infinitely differentiable function with bounded derivatives. Then $y_n$ can be 
represented by the series in powers of $\tau$:
\[
y_n = \sum_{m = 0}^{\infty} a_{m,n}\tau^m + O(\tau^k),
\]
where $k$ is the order of approximation.
\end{theorem}
This follows from the Taylor expansion of $G(y_n)$ at $y^*=a_{0,n}$:
\[
G(y_n) = G(y^*) + \sum_{m=1}^{\infty}\frac{G^{(m)}(y^*)}{m!}\left(y_n-y^*\right)^m
\]

Using this expansion we get $y_{n+1}$. Comparison of the coefficients gives the following
recurrent relations:
\begin{equation}\label{eq:recurrent_coefficients}
\begin{array}{l}
a_{0,n+1} = a_{0,n}, \\
a_{1,n+1} = a_{1,n} + G(y^*), \\
a_{2,n+1} = a_{2,n} + a_{1,n}G'(y^*), \\
a_{3,n+1} = a_{2,n} + a_{2,n}G'(y^*) + \frac{1}{2}a_{1,n}^2G''(y^*),
\end{array}
\end{equation}
and so on.

The expansion of $y_n$ can not be truncated after a fixed number of terms because the 
truncated value will be of greater order than the accuracy of the finite-difference method, such
a solution would diverge.

\end{comment}

Let us construct a procedure of redistribution of the values of the coefficients $a_i$ so that

1) the value of $y_n$ remains unchanged;

2) for any $i$: $|a_i|\leq1/\tau$.  

Consider the following example. Let $\tau = 0.1$ and in $n$-th layer $y_n$ is written as follows:
\[
y_n = 5\tau^0 + 2 \tau^1 + 6\tau^2 = 5.26.
\]
Let after applying the finite-difference formula we got $y_{n+1}$
\[
y_{n+1} = 5\tau^0 + 2 \tau^1 + 15\tau^2 = 5.35.
\]
The coefficient at $\tau^2$ is greater than $1/\tau = 10$ which leads to the eventual divergence. In order to 
avoid this we apply the digit shifting procedure.
It can be done in different ways, but following to how we do it in common positional numeral systems, we
take $10$ from $15$ and add $10\tau^2 = \tau$ to $2$ (at $\tau^1$):
\[
y_{n+1} = 5\tau^0 + 3 \tau^1 + 5\tau^2 = 5.35.
\]
As we see, the value is not changed, but the structure of the value is changed. 

The digit shifting procedure 
allows us to hold the values when applying the numerical scheme which would otherwise be 
truncated.

\begin{figure}[H]
\includegraphics[width=10cm]{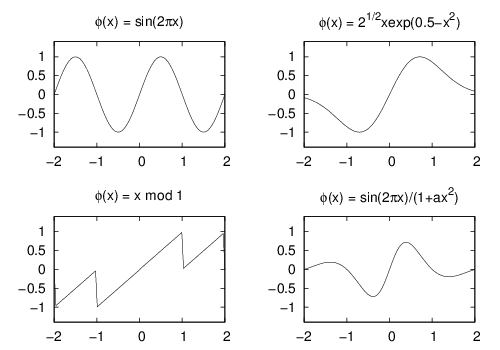} 
\caption{Shifting functions.}
\label{shifters}
\end{figure}

Now we will show how the shifting procedure can be constructed. To guarantee that series 
\eqref{eq:y_series} is convergent we will restrain the absolute values of $a_i$ by $\tau^{-1}$. 
Let there is a piecewise 
function $\phi(x)$ that maps $x$ to the segment $[-1, 1]$. Consider a function 
$\psi(x)=\tau^{-1}\varphi(\tau x)$  
We will refer to it as \emph{the shifting function}. In Fig. \ref{shifters} 
few possible choices of $\varphi(x)$ are shown.

Consider the expansion $y_n$ where the absolute values of the coefficients can be greater than $\tau^{-1}$:
\[
y_n = \tilde{a}_{0,n} + \tilde{a}_{1,n}\tau + \tilde{a}_{2,n}\tau^2 + ...
\]
Let $\phi(x)$ is the shifting function. It should be applied for each coefficient. 
Let us consider its application for $k$-th coefficient. First of all we must add to it the
value shifted from $k+1$-th coefficient, denote it as $\delta_{k+1,n}$. Let us add and subtract
the value $\psi(\tilde{a}_{k,n}+\delta_{n,k+1})$
\[
\psi(\tilde{a}_{k,n}+\delta_{k+1,n})+ \tilde{a}_{k,n}+\delta_{k+1,n}- 
\psi(\tilde{a}_{k,n}+\delta_{k+1,n}).
\]
Since $|\psi(x)|\leq 1/\tau$, we imply the following:
\[
a_{k,n} \equiv \psi(\tilde{a}_{k,n}+\delta_{k+1,n}),
\]
\[
\delta_{k,n}\equiv \tau\left(\tilde{a}_{k,n}+\delta_{k+1,n}- 
\psi(\tilde{a}_{k,n}+\delta_{k+1,n})\right).
\]

The shifting procedure guarantees that the absolute values of the coefficients 
of the power series are limited by
the value $\tau^{-1}$ which implies that the series may be truncated. 

\subsection{General solution representation}
Consider series \eqref{eq:y_series}. Now, using the digit shifting procedure,
we will construct the representation of the solution. Let us denote $a_{0,n} + a_{1,n}\tau$ as
$\alpha$ and from $a_{2,n}$ and on as $\xi$. 
Since all coefficients $a_{i,n}$ are limited by $\tau^{-1}$, we conclude that
\[
\xi < \tau  + \tau^2 + ... = O(\tau),
\]
thus $s$ is of order of $\tau$. Let us show that $f(y_n)$ can be approximated by 
$f(\alpha)$ with error $O(\tau)$. Consider the difference $f(y_n) - f(\alpha)$. 
Since $y_n \equiv \alpha + \xi$, let us expand $f(y_n)$ at $\alpha$
\[
f(y_n) = f(\alpha) + f'(\alpha)\xi + O(\xi^2).
\]
Since the derivatives of $f$ are bounded, we obtain that $f(y_n) = f(\alpha) + O(\xi)$. This
yields $f(y_n) = f(a_{0,n} + a_{1,n}\tau) + O(\tau)$.

Let us expand $f(a_{0,n} + a_{1,n}\tau)$ at $a_{0,n}$:
\begin{equation}\label{eq:fapprox}
f(a_{0,n} + a_{1,n}\tau) = f(a_{0,n}) + \sum_{m=1}^\infty\frac{f^{(m)}(a_{0,n})}{m!}(a_{1,n}\tau)^m.
\end{equation}

Consider the first order explicit finite-difference scheme for \eqref{eq:ode}:
\[
y_{n+1} = y_n + f(y_n).
\]
Combining \eqref{eq:y_series} and \eqref{eq:fapprox} we obtain:
\[
y_{n+1} = a_{0,n} + a_{1,n}\tau + a_{2,n}\tau^2 + 
\tau\sum_{m=0}^\infty\frac{f^{(m)}(a_{0,n})}{m!}(a_{1,n}\tau)^m.
\]
Then we get:
\begin{multline}\label{eq:yn_approx}
y_{n+1} = a_{0,n} + \left(a_{1,n} + f(a_{0,n})\right)\tau + 
\left(a_{2,n} + f'(a_{0,n})a_{1,n}\right)\tau^2 + \\ +
\sum_{m=2}^\infty\frac{f^{(m)}(a_{0,n})}{m!}a_{1,n}^m\tau^{m+1}.
\end{multline}

\begin{theorem}
If the right hand part of ODE is a polynomial of degree $r$, then for the first order explicit scheme
$p = r + 1$. 
\end{theorem}

A proof follows from the fact that the derivatives of the order of $r+1$ and greater equals zero for 
a polynomial of degree $r$. This yields that for any $m>r$ in \eqref{eq:yn_approx} the 
corresponding term equals zero, thus \eqref{eq:yn_approx} can be truncated at term with $\tau^{r+1}$.

\begin {comment}

Let us expand $f(y_n)$ at $y* = a_{0,n}$. The derivatives of order $r+1$ and greater
will be equial to zero, thus
\[
f(y_n) = f(a_{0,n}) + f(a_{0,n})\left(a_{1,n}\tau + \ldots + a_{p,n}\tau^p\right) +
\ldots + f(a_{0,n})\frac{\left(a_{1,n}\tau + \ldots + a_{p,n}\tau^p\right)^r}{r!}
\]
Each $i$-th term contains the value $a_{1,n}\tau$ in power $i$. The absolute values of the coefficients 
are bound by the value $\tau^{-1} - 1$ by the digit shifting function. This leads that in each term
the only value of order of 1 is $a_{1,n}\tau$ which can not be truncated. Since the power of the last 
term is $r$ and $f(y_n)$ is multiplied by $\tau$ in finite difference scheme we conclude that

Without loss of generality we can ommit in $G$ the terms with powers of $\tau$ less than $r$ and
let $\gamma_r = 1$. The senior power of $\tau$ is defined by the senior terms in 
\eqref{eq:ynsinglelayer}. 
Let $0<k\leq q$ is the number of such a term that
\[
\tau(a_k\tau^k)^r = c\tau^q.
\]
This term (among others) will provide the greatest possible power of $\tau$. Without loss of 
generality we let $c = 1$. Thus
\begin{equation}\label{eq:powerequation}
a_k^r\tau^{kr+1} = \tau^q.
\end{equation}
The coefficients $a$ is of the order of $\tau^{-1}$ so we let $a_k = \tau^{-1}$:
\[
\tau^{kr+1-r} = \tau^q.
\]
This yields that
\[
k = \frac{q + r - 1}{r}
\]
and from \eqref{eq:powerequation} we get that the senior power of $\tau$ will be
\[
p = kr + 1 = q + r.
\]
Which proves the theorem.
\end{proof}
\end{comment}

\subsection {Solution splitting}
Consider a problem \eqref{eq:ode}. Let there is a solution  
\begin{equation}\label{eq:ynsinglelayer1}
y_n = a_0 + a_1\tau + ... + a_p\tau^p.
\end{equation}
As we have already shown the convergence is achieved by applying the shifting functions to 
constrain the coefficients. 
But instead of applying the shifting function to each coefficient, we can split 
the segment of the series into two
segments:
\begin{equation}
y_n = a_0 + a_1\tau + ... + a_q\tau^q + a_{q+1}\tau^{q+1} + \ldots +  a_p\tau^p=
\end{equation}
\begin{equation}
= (a_0\tau^{-q} + a_1\tau^{-q+1} + ... + a_q)\tau^q + (a_{q+1}\tau^{q+1-p} + \ldots +  a_p)\tau^p.
\end{equation}
This allows to write the solution as the following sum:
\[
y_n = \alpha_n\tau^q + \beta_n\tau^p,
\]
where $q$ is the order of approximation and $p$ is the number of retained terms. These two parts are
linked together via a single value shifting.

This approach simplifies constructing the solution. 
For example, $q$ could be equal to $1$ and $p$ could be the number greater, say $2$,
in this case the first term, namely $\alpha_k\tau^q$ describes the regular part of the solution and the second term $\beta_n\tau^p$ 
would be the quasirandom part.

\subsection {The basic shifting function}
Let the shifting function is $\psi(x) = x mod \tau^{-1}$. Consider $a_{m,n}$:
\[
a_{m,n} = \left(\tilde{a}_{m,n} + \delta_{m+1, n}\right)mod\tau^{-1}.
\]
This expression conforms to the congruent random number generator formula:
\[
x_m=(bx_{m-1} + c)mod P, \quad 0\leq b,c < p.
\]

We assume that coefficients $a_{q+1}, ..., a_p$ are the uniformly distributed random integer 
numbers. Therefore
the solution can be written as a sum of two parts:

deterministic: $a_{0,n} + ... + a_{r,n}$,

random: $a_{q+1,n}, ..., a_{p,n}$.

These two parts are linked together by the value $\delta_{q+1}$. If the probabilities of the 
values of $\delta_{q+1}$ are known, then we can exclude the intermediate layers of independent variable 
on which the deterministic part remains unchanged.

\section {Linear approximation}

Let us see how the method of computer analogy can be used to construct linear 
analytical approximations which in certain cases can be used to get qualitative (or maybe in some cases the quantitative) information
about the solution.
Let $p=1$, thus the solution will contain only two terms. Denote them as $a_n$ and $b_n$, thus:
\[
y_n = a_{n} + b_{n}\tau.
\]
Using \eqref{eq:recurrent_coefficients} we get the value in the next layer 
(let $|G(a_{n})|<\tau^{-1}$):
\[
y_n = a_{n} + (b_{n} + G(a_{n}))\tau.
\]
Now we apply the digit shifting procedure:
\[
y_{n+1} = (a_{n} + \delta_n) + \phi(b_{n} + G(a_{n}))\tau.
\]
Let $\phi(x) = x mod \tau^{-1}$. We will use this shifting function by default. The choice of this  
function is based on the fact that
for $x < \tau^1$ it follows that $\phi(x) = x$. Then $y_{n+1}$ can be written as follows:
\[
y_{n+1} = (a_{n} + \delta_n) + \left((b_{n} + G(a_{n}))\tau^{-1}\right)\tau, 
\quad \delta_n=[(b_n+G(a_n))\tau].
\]
Since $|G(a_{n})|<\tau^{-1}$, then possible values for $\delta_n$ are $-1, 0, 1$.

Comparing the coefficients on sequential layers we get the following reccurent relations:
\[
a_{n+1} = a_{n} + \delta_n, \quad b_{n+1} = (b_{n} + G(a_{n}))mod\tau^{-1}.
\]
We see that $a$ can be changed only when a non-zero digit shifting occures.  
Suppose that for several sequential layers $\delta$ equals to zero, so there are no values shifted to
coefficient $a$. This implyes that for these layers $G(a_{n})$ is constant and $b$ grows
linearly. When the digit shifting takes place the coefficient $a_0$ changes, that leads to changing
of angular coefficient of $b$ which is $G(a)$. 

When $b$ reaches $\pm\tau^{-1}$, $\delta_n$ will become $-1$ or $1$, this will change 
$a$ and consequently the angular coefficient $G$.

Let us now obtain the formulae for constructing the linear approximation. Let 
\[
y_0 = a_{0} + b_{0}\tau
\]
corresponds to some moment $t_0$. Suppose $\delta_0, \delta_1, ..., \delta_{M-1}$ equal to zero, so
there are no shiftings in the first $M$ layers, and there is a shift in $M-th$ layer. Then we imply
that $a_{m} = a_{0}$ for any $m$ that $0\leq m \leq M-1$. 

This means
that due to the chosen shifting function for any $m$ such that $0\leq m \leq M-1$ $\phi(x) = x$ because
the absolute value of the argument is less than $\tau^{-1}$. This yields the reccurent 
relation for $b$ can be written as follows:
\[
b_{m+1} = (b_m + G(a^*))mod\tau^{-1} = b_m+G(a_0).
\]
This relation can be written explicitly:
\[
b_{m} = b_0 + mG(a_0).
\]
By substituting $b_m$ to $y_m$ we get:
\[
y_m = a_m + b_m\tau = a_0 + b_0\tau+m\tau G(a_0).
\]
Notice that $m\tau$ is the value of time $t_m$, thus:
\[
y_m = y_0 + (t_m - t_0)G(a_0).
\]
This relation holds true for any value of $t$ such that $t_0 \leq t \leq t_M$, thus
\[
y = y_0 + (t - t_0)G(a_0).
\]
Now we will show the value $M$ can be obtained. Denote $G \equiv G(a_0)$ and consider two cases: 
$G\ge 0$ and $G < 0$. The first case implies that $b$ is linearly grows, the second one implies that
$b$ is linearly descends. Let $G \ge 0$ then the digit shifting occure when $b_M\ge \tau^{-1}$: 
\[
b_0 + MG \ge \tau^{-1}.
\]
Thus 
\[
M \ge \frac{1-b_0\tau}{G\tau}.
\]
For $M$ we will use the smallest integer value of the right hand part, so \mbox{$M \equiv \lceil\frac{1-b_0\tau}{G\tau}\rceil$}.
By analogy we get the following value for $M$ in a case $G<0$: $M \equiv \lceil-\frac{1+b_0\tau}{G\tau}\rceil$.
And now we need to show how to get the new values of $a_0$ and $b_0$ after the shifting of the digit.
Denote them as $a_0^*$ and $b_0^*$. If digit shifting is triggered, $b$ is changed according
to the shifting function, so:
\[
b_0^* = b_M mod \tau^{-1}.
\]
And $a$ is changed on $1$ or $-1$ according to $G$:
\[
a_0^* = a_0 + sign G.
\]

\section {Examples}

\begin{example}
Consider the following Cauchy problem (see [2] for details):
\[
\left\{
\begin{array} {l}
du/dt = v^2-u^2, \quad u(0) = 1 \\
dv/dt = u^2-2v, \quad v(0) = 0
\end{array}
\right.
\]
We will use the explicit first-order scheme and represent $u_n$ and $v_n$ as the segments of the
series in powers of $\tau$. 
\[
u_{n+1} = u_n + v_n^2\tau - u_n^2\tau, \quad
v_{n+1} = v_n + u_n^2\tau - 2v_n\tau.
\]

The proposed method allows us to obtain the explicit expressions:

\[
u_a = 1 - a\tau, \quad t_a = \tau\sum_{m=0}^{a-1}(1-m\tau)^{-2},\quad
v_a = \tau\sum_{m=1}^{a-1}\prod_{k=m+1}^{a-1}\left(1-\frac{2\tau}{(1-k\tau)^2}\right).
\]

The solutions comparison is shown on Fig. \ref{figCarlem}.

For small values of $t$ we can neglect the value $a^2\tau^2$ (see Fig. \ref{figSmallApprox}):
\[
u_a = 1-a\tau, \quad v_a = (1-2a\tau)\sum_{m=2}^a \frac{\tau}{1-2m\tau},
\quad t_a=\sum_{m=0}^{a-1} \frac{\tau}{1-2m\tau}.
\]
\begin{figure}[H]
\includegraphics[width=10cm]{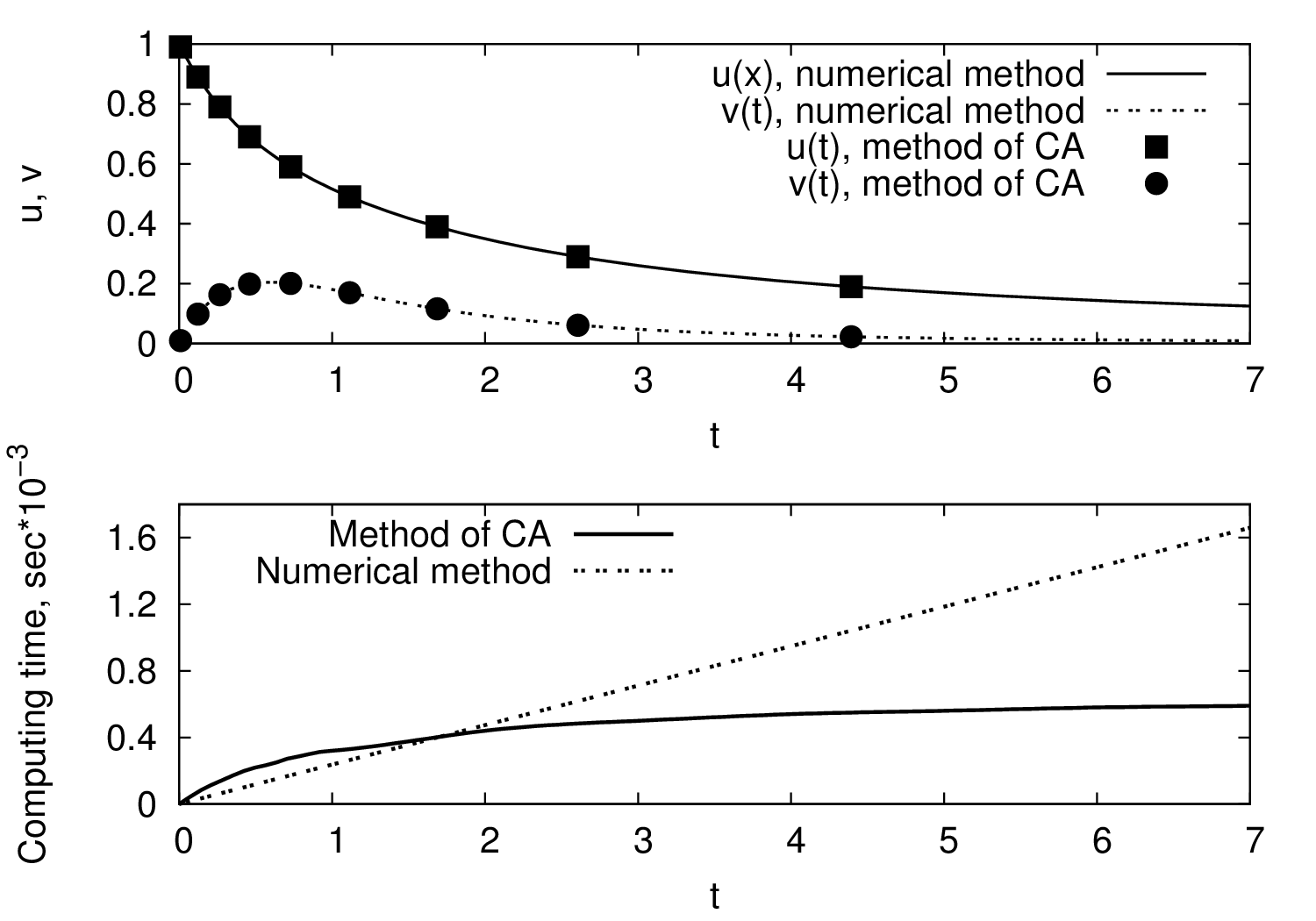} 
\caption {Solutions comparison}
\label{figCarlem}
\end{figure}

\begin{figure}[H]
\includegraphics[width=10cm]{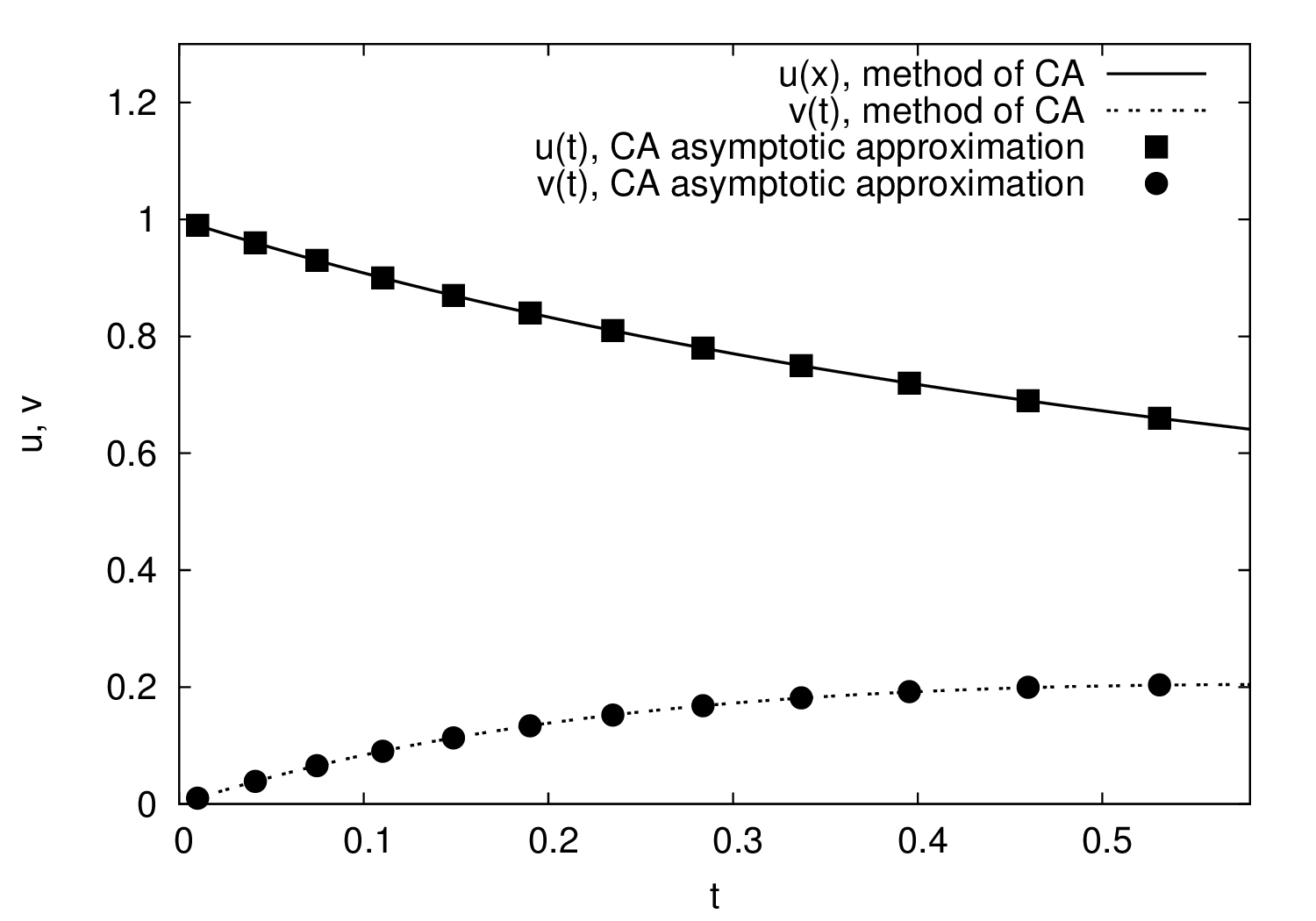} 
\caption {Approximation for small $t$}
\label{figSmallApprox}
\end{figure}

\end{example}

\begin{example}
Consider the Van der Pol oscillator which is a second order differential equation:
\[
\ddot u-(\lambda-u^2)\dot{u} + u = 0.
\]
For this problem there is no analytical solution. Let $\lambda = 1$ and rewrite 
the equation as a system of two differential equations:
\[
\left\{
\begin{array}{l}
\dot u = v, \\
\dot v = \left(1-u^2\right)v - u.
\end {array}
\right.
\]
Let $u(0) = 0$ and $v(0) = 1$. 
The first order finite difference solution will be as follows:
\[
\left\{
\begin{array}{l l}
u_{n+1} = u_n + v_n\tau, & \quad u_0 = 0, \\
v_{n+1} = v_n + v_n\tau - u_n^2v_n\tau-u_n\tau, & \quad v_0 = 1.
\end{array}
\right.
\]

We will represent the solution as a segment of the series in powers of $\tau$. 
Theorem 1 implies that for $u$ and $v$ it is sufficient to retain terms until $\tau^2$ and 
$\tau^4$ respectively. We will retain only the terms up to $\tau^1$ for the both functions:
\[
\left\{
\begin{array}{l}
u_n = a_{0,n} + a_{1,n}\tau, \\
v_n = b_{0,n} + b_{1,n}\tau.
\end{array}
\right.
\]
Taking into account the mentioned notices, this solution will not converge, but it is still useful
because it allows us to construct linear a approximation which reproduces the behaviour of 
the solution.
Using the finite-difference scheme and the digit shifting procedure we obtain:
\[
\left\{
\begin{array}{l}
u_{n+1} = (a_{0,n} + \delta_n) + \left((a_{1,n}+b_{0,n})mod\tau^{-1}\right), \\
v_{n+1} = (b_{0,n} + \omega_n) + 
\left((b_{1,n} + b_{0,n} - a_{0,n} - a_{0,n}^2b_{0,n})mod\tau^{-1}\right).
\end{array}
\right.
\]
Here the shifted values are given by the following expressions:
\[
\delta_n = [(a_{1,n} + b_{0,n})\tau], 
\quad \omega_n = [(b_{1,n} + b_{0,n} - a_{0,n} - a_{0,n}^2b_{0,n})\tau].
\]

Let us see how a linear approximation can be constructed. From the initial condition we 
imply that 
\[
a_{0,0} = 0, \quad a_{1,0} = 0, \quad b_{0,0} = 1, \quad b_{1,0} = 0.
\]

\emph{Layer 0}

From the initial values we obtain $G_1$ and $G_2$:
\[
G_1(a_{0,0}, b_{0,0}) = b_{0,0} = 1, \quad G_2(a_{0,0}, b_{0,0}) = (1-a_{0,0}^2)b_{0,0} - a_{0,0} = 1.
\]
Now we calculate $M_1$ and $M_2$:
\[
M_1 = \frac{1-a_{1,0}\tau}{G_1\tau} = \frac{1}{\tau}, \quad
M_2 = \frac{1-b_{1,0}\tau}{G_2\tau} = \frac{1}{\tau}.
\]
A value $M$ will be the minimum of values $M_1$ and $M_2$, namely:
\[
M = min(M_1, M_2) = 1/\tau.
\]
A value $M$ defines $t_M$:
\[
t_M = t_0 + M\tau = 0 + 1 = 1.
\]
Now we can write the linear approximation in a segment $[0, 1]$:
\[
\left\{
\begin{array}{l}
u(t) = u_0 + G_1(t-t_0) = 0 + 1(t-0) = t, \\
v(t) = v_0 + G_2(t-t_0) = 1 + 1(t-0) = 1 + t. 
\end{array}
\right.
\]
To get a solution in the next layer we must compute the initial values for it. Since $G_1>0$ and $G_2>0$ the
new initial values for $a_0$ and $b_0$ will be $a_{0,0}^*$ and $b_{0,0}^*$:  
\[
a_{0,0}^* = a_{0,0} + 1 = 1, \quad b_{0,0}^* = b_{0,0} + 1 = 2.
\]
The initial values for $a_1$ and $b_1$ are found from $a_{1,M}$ and $b_{1,M}$:
\[
a_{1,0}^* = a{1,M}mod\tau^{-1} = 0, \quad b_{1,0}^* = b{1,M}mod\tau^{-1} = 0.
\]

\emph{Layer 1}

By analogy with Layer 0 we compute $G_1$ and $G_2$:
\[
G_1 = 2, \quad G_2 = -1.
\]
Find $M_1$ and $M_2$:
\[
M_1 = 1/\tau \quad M_2 = 1/2\tau.
\]
This yields $M = 1/2\tau$ and $t_M = 1 + 0.5 = 1.5$ and we construct
the following approximation in a segment $[1, 1.5]$:
\[
\left\{
\begin{array}{l}
u(t) = u_0 + G_1(t-t_0) = 1 + 2(t-1), \\
v(t) = v_0 + G_2(t-t_0) = 2 - (t-1). 
\end{array}
\right.
\]
The initial conditions for the next layer are as follows:
\[
a_{0,0}^* = a_{0,0} + 1 = 2, \quad b_{0,0}^* = b_{0,0} = 2.
\]
A coefficient $b_{0,0}$ is not changed because $M_2 < M_1$, so no shifting in $v$ occures. We note that this linear solution does not depend on a magnitude of the step $\tau$.

\emph{Layer 2}

Values of function $G$ are:
\[
G_1 = 2, \quad G_2 = -8.
\]
And we obtain $M_1$ $M_2$:
\[
M_1 = 1/4\tau, \quad M_2 = 1/8\tau.
\]
We get the following approximation in a segment $[1.5, 1.625]$:
\[
\left\{
\begin{array}{l}
u(t) = u_0 + G_1(t-t_0) = 2 + 2(t-1.5), \\
v(t) = v_0 + G_2(t-t_0) = 0.5 - 8(t-1.5). 
\end{array}
\right.
\]
The explicit solutions can be written in the following table: 
\begin{center}
\begin{tabular}{|l|l|l|l|l|l|}
\hline
$t$ & $u(t)$ & $v(t)$ & $t$ & $u(t)$ & $v(t)$ \\
\hline
$[0, 1]$ & $t$ & $1+t$ &
$[6.25, 6.4]$  & $-2t+10.5$ &  $ 8t-52.18  $ \\ 

$[1, 1.5]$ & $2t-1$ & $-t + 3$ &
$[6.4, 6.6]$ &   $-1t+4.1$ &   $5t-32.98  $ \\ 

$[1.5, 1.57]$ & $t-1$ & $-8t+13.5$ &
$[6.6, 7.09]$ &  $-2.5$ &   $2t-13.18  $ \\ 

$[1.57, 1.76]$ & $t+0.57$ &   $-5t+8.79  $ &
$[7.09, 8.09]$ & $t-9.59$ & $  -t+8.09  $ \\ 

$[1.76, 2.26]$ & $2.33$ &   $-2t+3.51  $ &
$[8.09, 8.59]$ & $-1.5$ & $  2t-16.18  $ \\ 

$[2.26, 3.27]$ & $-t+4.59$ &  $ t-3.27  $ &
$[8.59, 9.09]$ & $t-10.09$ & $  -t+9.59$ \\   

$[3.27, 3.77]$ & $1.32$ &   $-2t+6.54  $ &
$[9.09, 10.09]$ &$t-10.09$ & $  t-8.59$ \\   

$[3.77, 4.09]$ & $-t+5.09$ &  $ t-4.77  $ &
$[10.09, 10.59]$&$t-10.09$ & $  t-8.59  $ \\ 

$[4.09, 5.09]$ & $-t+5.09$ &  $ -t+3.41  $ &
$[10.59, 10.84]$&$2t-20.68$ & $  2t-19.18$ \\   

$[5.09, 5.41]$ & $-t+5.09$ &  $ -t+3.41  $ &
$[10.84, 11.34]$&$2t-20.68$ & $  -t+13.34$ \\   

$[5.41, 5.75]$ & $-2t+10.5$ &  $ -2t+8.82  $ & 
$[11.34, 11.47]$&$2t-20.68$ & $  -8t+92.72$ \\   

$[5.75, 6.25]$ & $-2t+10.5$ &  $ t-8.43  $ &
$[11.47, 11.67]$&$t-9.21$ &  $ -5t+58.31 $ \\
\hline
\end{tabular}

\textbf{Table 1.} Piecwise linear approximation for Van der Pol equation.
\end{center}
Figures \ref{figVanDerPolLinearSingle} and \ref{figVanDerPolLinearMulti} demonstrate that this solution converges to the cycle which roughly corresponds to the 
exact solution.
\begin{figure}[H]
\includegraphics[width=10cm]{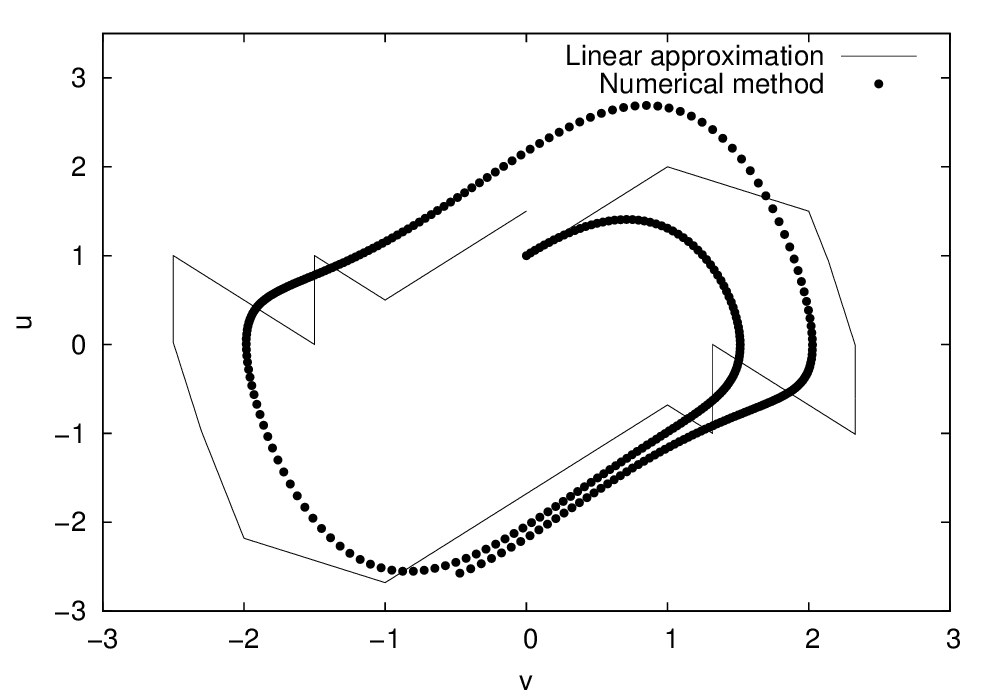} 
\caption {Linear approximation, a single cycle.}
\label {figVanDerPolLinearSingle}
\end{figure}
Thus we obtain an analytical solution of the first approximation which provides a qualitative information and  
reflects some quantitative features.
\begin{figure}[H]
\includegraphics[width=10cm]{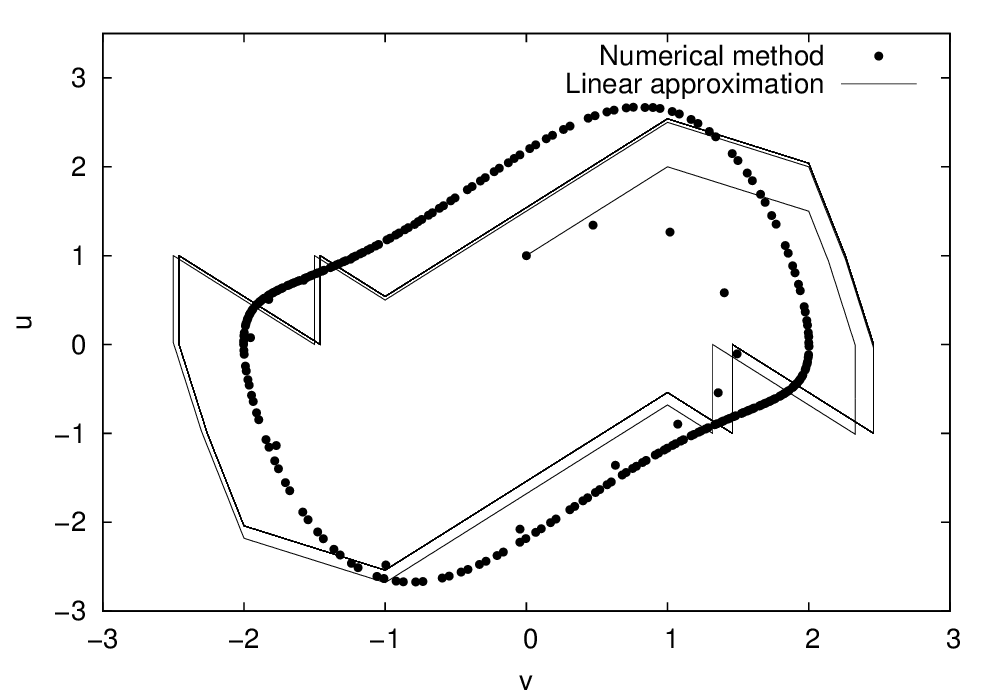} 
\caption {Linear approximation, multiple cycles.}
\label {figVanDerPolLinearMulti}
\end{figure}
\end{example}

\begin{example}
The Lorenz system is given by the following system of nonlinear differential equations:
\[
\left\{
\begin{array} {l}
dx/dt = \sigma(y-x), \quad x(0) = 3 \\
dy/dt = x(r-z) - y, \quad y(0) = 2 \\
dz/dt = xy-vz, \quad z(0) = 15.
\end{array}
\right.
\]
Let us use the following parameters (which lead to solution with a strange attractor): $\sigma = 3, \; r=15, \; v = 1$. We represent $x_n, y_n, z_n$ according to 2.4:
\[
\begin{array} {l}
x_{n} = a_{1,n}\tau + a_{2,n}\tau^3,
y_{n} = b_{1,n}\tau + b_{2,n}\tau^3, 
z_{n} = c_{1,n}\tau + c_{2,n}\tau^3.
\end{array}
\]

The finite-difference method and the shifting procedure give:
\[
a_{1,n+1} = a_{1,n} + \delta_n, \quad 
a_{2,n+1} = \left(a_{2,n} + \frac{\sigma}{\tau}(b_{1,n}-a_{1,n})\right)mod\tau^{-2}, 
\]
\[
b_{1,n+1} = b_{1,n} + \omega_n,  \quad 
b_{2,n+1} = \left(b_{2,n} + r\frac{a_{1,n}}{\tau}-a_{1,n}c_{0,n}-\frac{b_{1,n}}{\tau}\right)mod\tau^{-2},
\]
\[
c_{1,n+1} = c_{1,n} + \eta_n, \quad 
c_{2,n+1} = \left(c_{2,n} + a_{1,n}b_{1,n}-v\frac{c_{1,n}}{\tau}\right)mod\tau^{-2} 
\]
\[
\delta_n = [\tau\sigma(b_{1,n}-a_{1,n})], \quad \omega_n = [\tau\left(ra_{1,n}-a_{1,n}c_{0,n}\tau-b_{1,n}\right)],
\]
\[
\eta_n = [\tau\left(a_{1,n}b_{1,n}\tau - vc_{1,n}\right)].
\]
Simulations confirm that this solution converges to the numerical solution.

Function $x(t), y(t), z(t)$ is shown in Fig. \ref{figLorenzTxyz}. Simulations using method of CA shows the Lorenz attractor behaviour (Fig. \ref{figSimulation1}).
\begin{figure}[H]
\includegraphics[width=13cm]{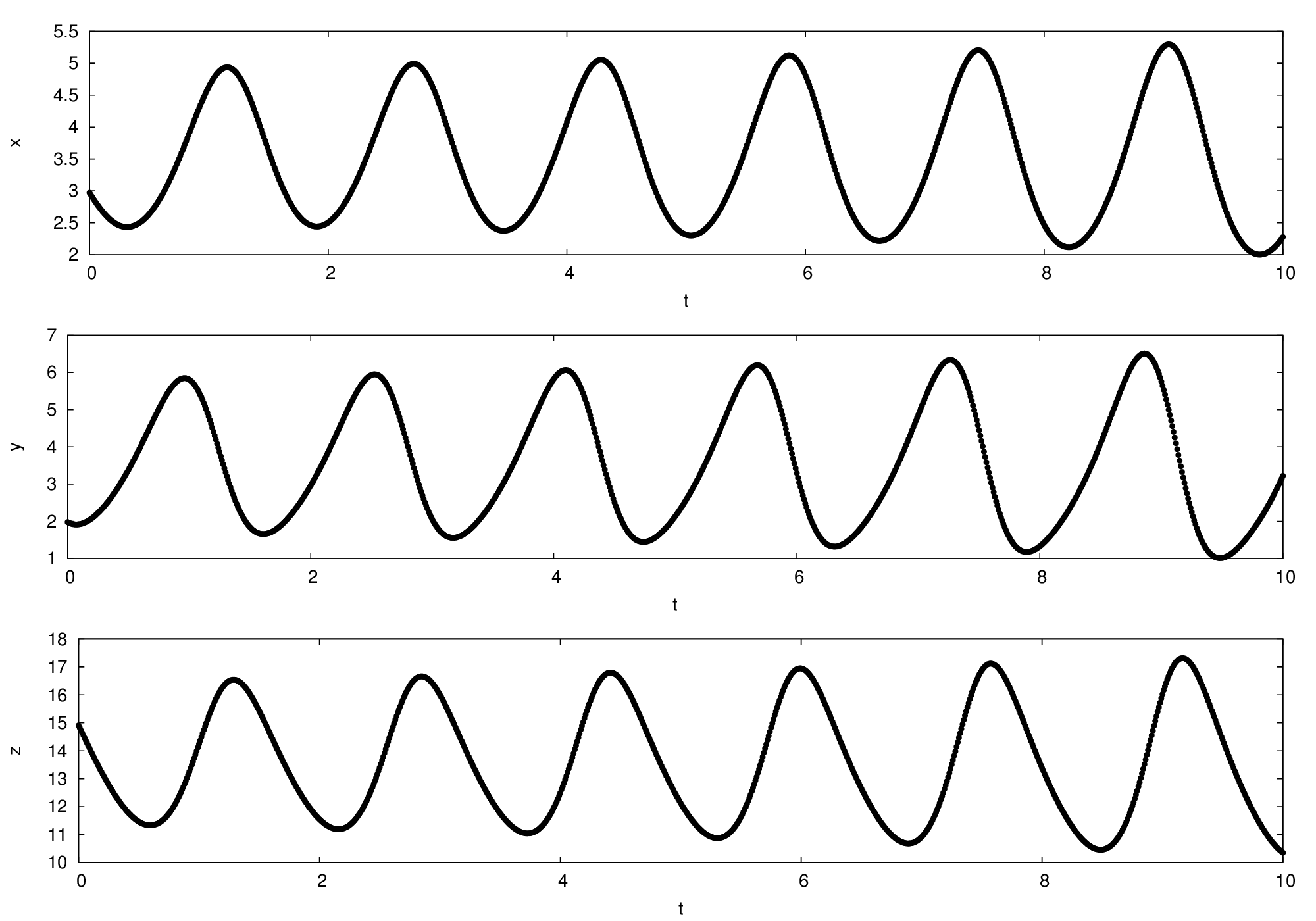} 
\caption {$x(t), y(t), z(t)$}
\label{figLorenzTxyz}
\end{figure}
\noindent
\begin{figure}[H]
\includegraphics[width=13cm]{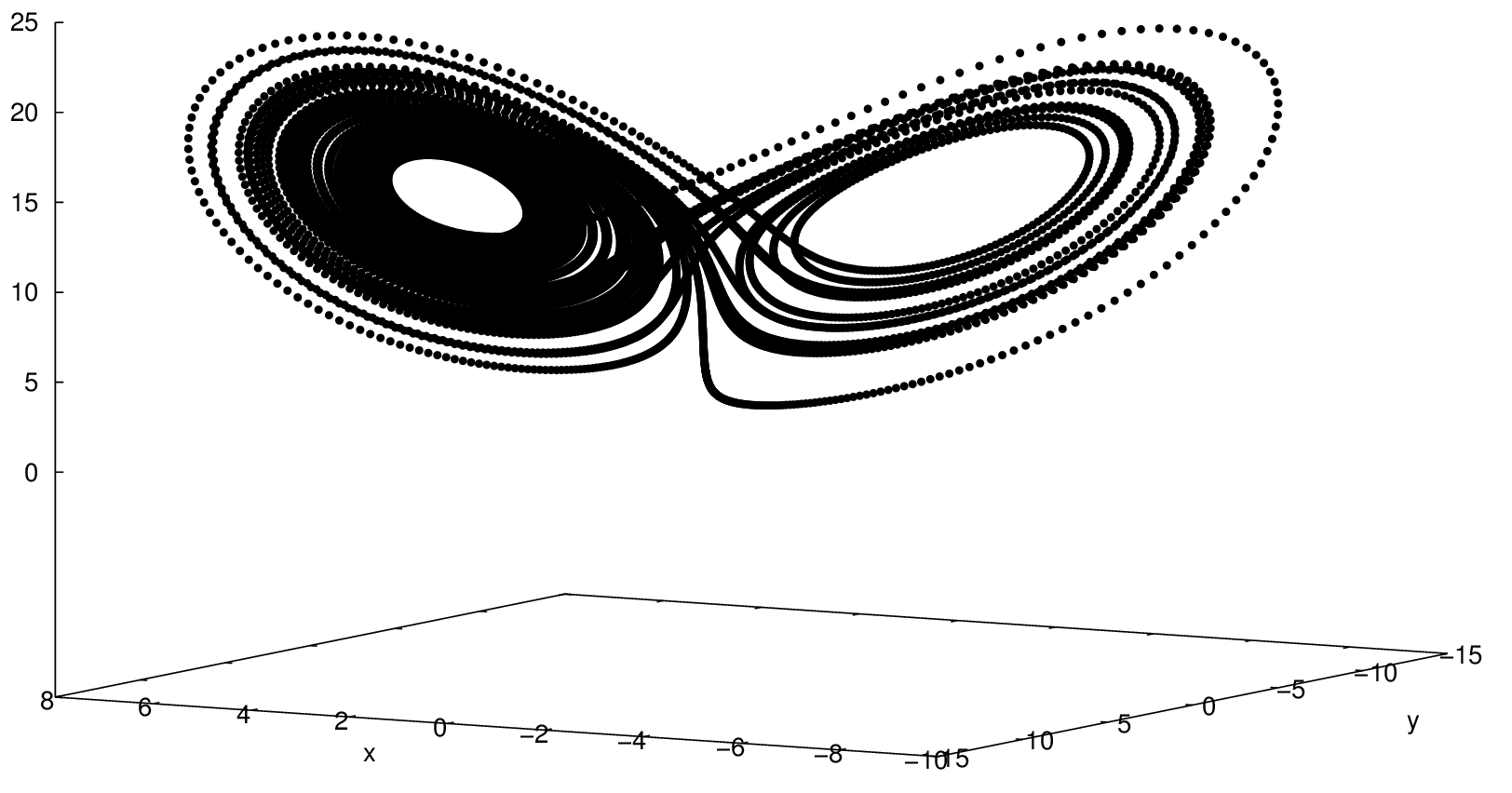} 
\caption {Simulation of the solution.}
\label{figSimulation1}
\end{figure}

Let us construct the linear approximation. We will retain
only the linear terms:
\[
\begin{array} {l}
x_{n+1} = (a_{0,n}+\delta_{n}) + 
\left((a_{1,n}+\sigma(b_{0,n}-a_{0,n}))mod\tau^{-1}\right)\tau
\\
y_{n+1} = (b_{0,n}+\omega_{n}) + 
\left((b_{1,n}+a_{0,n}r-a_{0,n}c_{0,n}-b_{0,n})mod\tau^{-1}\right)\tau
\\
z_{n+1} = (c_{0,n} + \eta_{n}) +  
\left((c_{1,n}+a_{0,n}b_{0,n}-vc_{0,n})mod\tau^{-1}\right)\tau
\end{array}
\]
Using the formulae from section 3 we get the following approximation
\begin{center}
\begin{tabular}{|c|c|c|c|c|c|c|c|}
\hline
t & x & y & z & t & x & y & z \\
\hline
$[0, 0.12]$ & $-3t+3$ & $-2t+2$ & $-9t+15$ & 
$[0.68, 0.71]$ & $3t+0.39$ & $5t+0.39$ & $-4t+12.71$ \\   

$[0.12, 0.24]$ & $-3t+3$ & $-2t+2$ & $-8t+14.88$ & 
$[0.71, 0.79]$ & $3t+0.39$ & $7t-0.97$ & $-2t+11.29$ \\  

$[0.24, 0.34]$ & $-3t+3$ & $t+1.64$ & $-7t+14.64$ &
$[0.79, 0.84]$ & $6t-1.74$ & $6t-0.26$ & $2t+8.13$ \\   

$[0.34, 0.37]$ & $-3t+3$ & $4t+0.92$ & $-9t+15.32$ &
$[0.84, 0.94]$ & $3t+0.63$ & $11t-4.21$ & $5t+5.61$ \\   

$[0.37, 0.5]$ & $1.98$ & $2t+1.6$ & $-8t+14.95$ &
$[0.94, 0.97]$ & $6t-1.89$ & $10t-3.37$ & $8t+2.79$ \\   

$[0.5, 0.53]$ & $1.98$ &  $4t+0.86$ & $-7t+14.45$ &   
$[0.97, 1.01]$ & $9t-4.71$ & $9t-2.43$ & $14t-3.03$  \\ 

$[0.53, 0.68]$ & $1.98$ &  $6t-0.14$ & $-5t+13.39$ & & & & \\   

\hline
\end{tabular}
\end{center}
This solution doesn't converge, but it still preserves some information about the exact solution. We
can see that there is still periodic behavior (Fig. \ref{figLorenzTxyzApprox}) which corresponds to the 
		converging solution.
\begin{figure}[H]
\includegraphics[width=13cm]{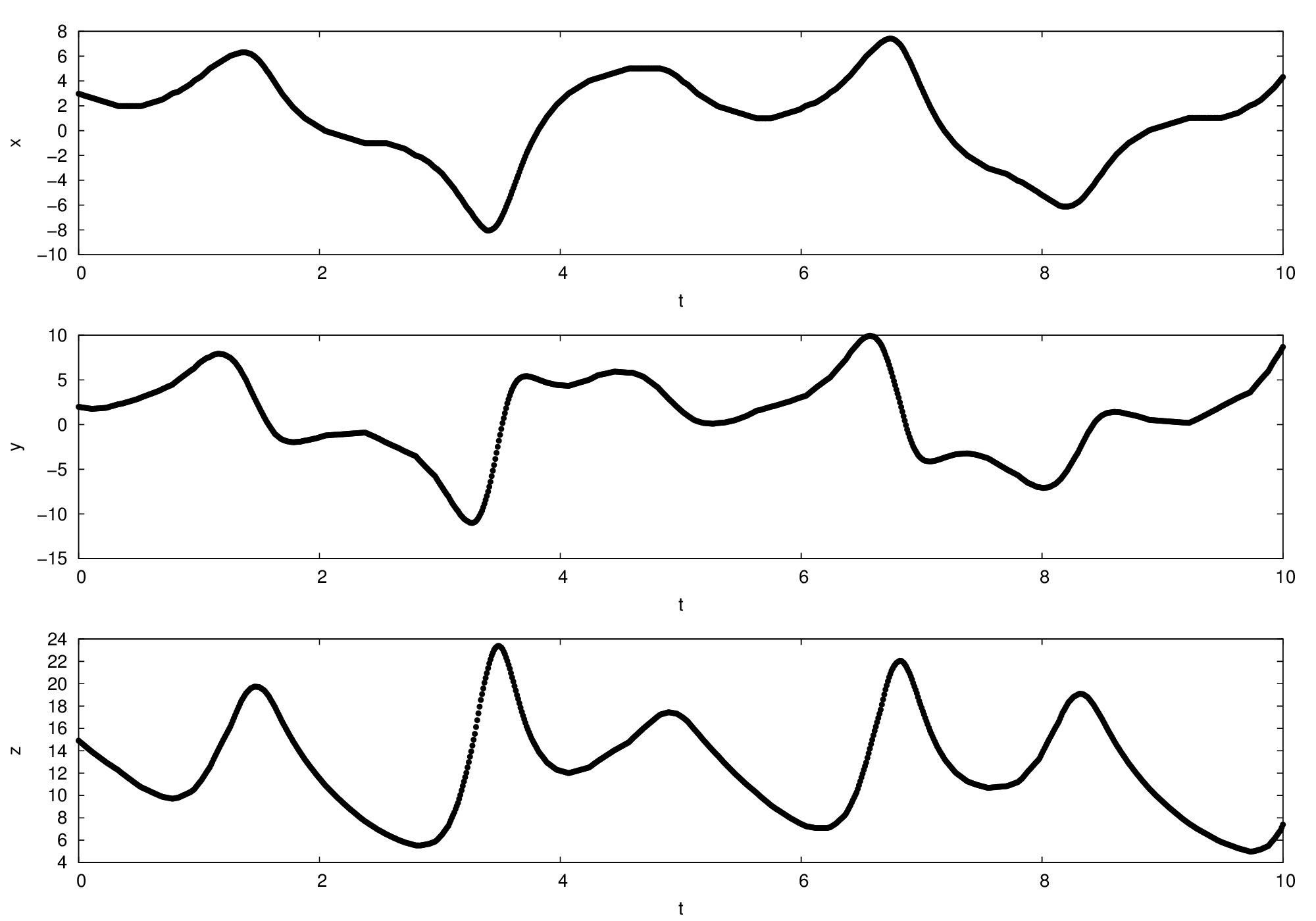} 
\caption{$x(t), y(t), z(t)$, linear approximation.}
\label{figLorenzTxyzApprox}
\end{figure}

Linear approximation of the Lorenz attractor is shown in Fig. \ref{figLorenzLinearApprox}.
\begin{figure}[H]
\includegraphics[width=13cm]{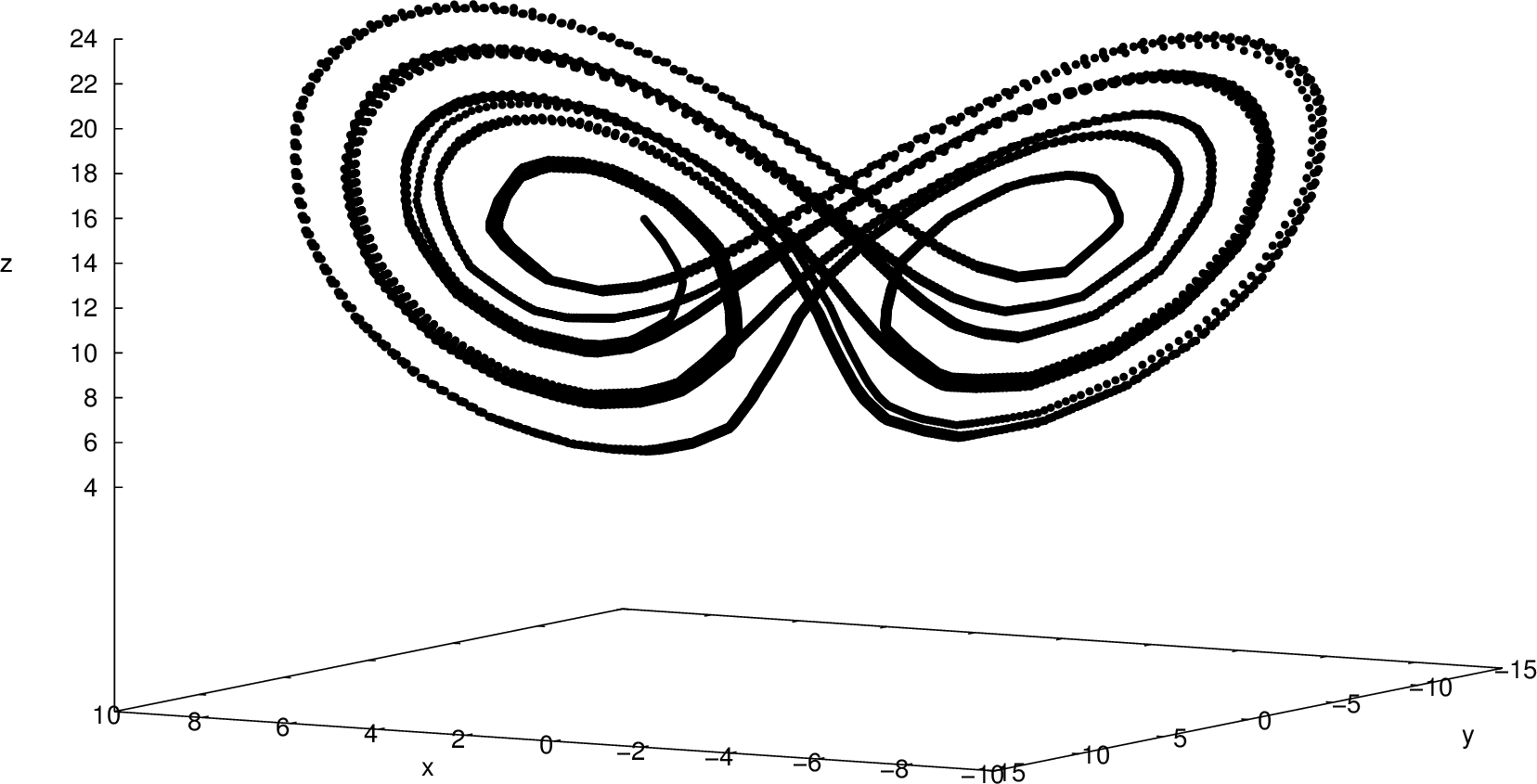} 
\caption {Linear approximation of the solution.}
\label {figLorenzLinearApprox}
\end{figure}
\end{example}

\section {Conclusions}
A theoretical model of the computer working with numbers was presented. We analyzed the possibility of formalizing 
the computer operations, and proposed a representation of a solution in the form of a segment 
of the series in the powers of the step of the independent variable. This technique involves 
analytical work for obtaining the probabilities. However, the main steps can be listed as 
follows: to choose a convergent finite-difference scheme, to find the number of terms to be 
retained, and to obtain the expressions for the main coefficients in the representation of the 
solution. The proposed method of the computer analogy allows us to represent the solution of 
the problem as a convergent series that can be analyzed.

%
%

\end{document}